\documentclass[12pt]{article}

\usepackage{amsmath,amsthm}
\usepackage{amssymb}
\usepackage{amsfonts}
\usepackage{amscd}

\newtheorem{thm}{Theorem}[section]
\newtheorem{theorem}[thm]{Theorem}

\newtheorem{corollary}[thm]{Corollary}
\newtheorem{lemma}[thm]{Lemma}

\theoremstyle{definition}

\newtheorem{definition}[thm]{Definition}

\newtheorem{examples}[thm]{Examples}

\newtheorem{observation}[thm]{Observation}

\newcommand\al{{\alpha}}

\newcommand\la{{\lambda}}
\newcommand\s{{\sigma}}

\newcommand\ep{{\epsilon}}

\newcommand\Z{{\bold Z}}
\newcommand\C{{\bold C}}
\newcommand\bX{{\bold X}}
\newcommand\bg{{\bold g}}
\newcommand\rr{{r}}
\newcommand\GL{{\bf GL}}
\newcommand\bV{{\bf V}}
\newcommand\bT{{\bf T}}

\newcommand\bG{{\bf G}}

\newcommand\fg{{\mathfrak g}}
\newcommand\fgl{{\mathfrak{gl}}}
\newcommand\bgl{{\bold{gl}}}

\newcommand\rGL{{\hbox{GL}}}
\newcommand\bGL{{\bold{GL}}}

\newcommand\Hom{{\hbox{Hom}}}
\newcommand\End{{\mathrm{End}}}
\newcommand\rExt{{\mathrm{Ext}}}
\newcommand\bEnd{{\bold{End}}}
\newcommand\defi{{\text{def}}}

\newcommand\re{{\hbox{e}}}
\newcommand\bfe{{\hbox{\bf e}}}

\newcommand\bE{{\hbox{\bf E}}}
\newcommand\bH{{\hbox{\bf H}}}
\newcommand\Sym{{\hbox{Sym}}}

\newcommand\sets{{ \hbox{(sets)} }}
\newcommand\groups{{ \hbox{(groups)} }}

\newcommand\salg{{(\mathrm{salg})}}

\newcommand\lra{{\longrightarrow}}

\newcommand\Proof{{\bf Proof}}

\begin{document}

\vskip 1cm 
\centerline{\LARGE \bf Tensor Representations of the }
\bigskip
\centerline{\LARGE \bf General Linear Super Group}
\bigskip

\centerline{R. Fioresi {\footnote{Investigation supported by
the University of Bologna, funds for selected research topics.}}}

\bigskip

\centerline{\it Dipartimento di Matematica, Universit\`a
di Bologna }
\centerline{\it Piazza di Porta S. Donato, 5. 40127 Bologna. Italy.}
\centerline{{\footnotesize e-mail: fioresi@dm.unibo.it}}

\bigskip

\begin{abstract}
We show a correspondence between tensor representations of the
super general linear group $GL(m|n)$ and
tensor representations of the general linear
superalgebra $\fgl(m|n)$ constructed by Berele and Regev
in \cite{br}.
\end{abstract}

\section{Introduction}

\medskip 

Supersymmetry is an important mathematical tool in physics
that enables to treat on equal grounds the two types
of elementary particles: bosons and fermions, whose states are
described respectively by commuting and anticommuting functions.
It is fundamental to seek a unified treatment for these particles
since they do transform into each other.  Hence a symmetry
that keeps one type separated from the other is not acceptable.
For this reason the symmetries of elementary particles 
must be described not by 
groups, but by {\it supergroups}, which are a natural generalization
of groups in the $\Z_2$ graded or {\it super} setting.

\medskip

The theory of representations of supergroups has a particular
importance since it is attached to the problem of the classification
of elementary particles. For a more detailed historical
and physical introduction to supersymmetry see the
beautiful treatment in \cite{vsv} 1.7, 1.8.

\medskip

As in the classical theory, in order to
understand the representations of a supergroup, one must first study
the representations of its Lie superalgebra.
The representation theory of the general linear superalgebra
$\fgl(m|n)$ has been the object of study of many people.

\medskip

In \cite{br} Berele and Regeev provide a full account of a class of
irreducible representations of $\fgl(m|n)$ that turns 
out to be linked to certain Young tableaux
called {\it semistandard or superstandard tableax}. The same
result appeares also in \cite{dj} by Dondi and Jarvis 
in a slightly different setting. Dondi and Jarvis
in fact introduce the notion of {\it super permutation} and
use this definition to motivate the semistandard Young
tableaux used for the description of the irreducible
representations of the general linear superalgebra.

\medskip
 
The results by Berele and Regeev were later generalized
and deepened by Brini, Regonati and Teolis in \cite{brt}.
In their important work, they develop a unified theory
that treats simoultaneosly the super and the classical
case, through the powerful method of {\it virtual variables}.

\medskip
 
Another account of this subject is found in
\cite{se}. Sergeev establishes a correspondence between
a class of irreducible tensor representations of $\fgl(m|n)$ and
the irreducible representations of a certain finite group, 
though different from the permutation group used both in 
\cite{br} and \cite{dj}.

\medskip

It is important to remark at this point that the theory
of representations of superalgebras and of supergroups
has dramatic differences with respect to the classical
theory. As we will see, not all representations of 
the super general linear group and its Lie superalgebra are found
as tensor representations. 
Moreover not all representations
are completely reducible over $\C$.


\medskip

In this paper we want 
to understand how representations
of the Lie superalgebra $\fgl(m|n)$ can be naturally
associated to the representation of the corresponding group
$GL(m|n)$. Though this fact is stated in the physicists works
as for example \cite{dj}, \cite{bbb}, it is never satisfactorily
worked out. Using \cite{br, dj} we are
able then to obtain a full classification of the irreducible tensor
representations of the general linear supergroup coming from
the natural diagonal action.
We will do this using the approach suggested by Deligne and Morgan in
\cite{dm}: using {\sl the functor of points}.

\medskip

In fact while in non commutative geometry in general
the geometric object is lost and the
only informations are retrieved through various algebras, like
the $C^{\infty}$ functions or the algebraic functions, naturally
associated to it; in supergeometry, using
the functor of points approach, one is able to recover
the geometric intuition, which otherwise would be lost.

\medskip

This paper is organized as follows.

\medskip

In section \ref{definitions} we review some of the basic
definitions of supergeometry. Since we will adopt the functorial
language we relate  our definitions to the other definitions   
appearing in the literature.

\medskip

In section \ref{berele} we recall briefly the results obtained
indipendently by Berele, Regev and  Dondi, Jarvis.
These results establish a correspondence between tensor 
representations of the permutation group and tensor
representations of the superalgebra $\fgl(m|n)$. Moreover we show
that the tensor representations of the Lie superalgebra $\fgl(m|n)$ 
do not exhaust all polynomial representations of $\fgl(m|n)$.

\medskip

Finally in section \ref{supergroup} we 
discuss tensor representations of the general linear supergroup
associated to the representations of $\fgl(m|n)$ described in
\S 3.

\medskip

{\bf Acknoledgements}. We wish to thank Prof. V. S. Varadarajan,
Prof. A. Brini, Prof. F. Regonati and Prof. I. Dimitrov
for helpful comments.

\bigskip

\section{Basic definitions} \label{definitions}

\medskip

Let $k$ be an algebraically closed field of characteristic 0. 

\medskip

All algebras have to be intended over $k$.

\medskip

A {\it superalgebra} $A$ is a $\Z_2$-graded algebra,
$A=A_0 \oplus A_1$, $p(x)$ will
denote the parity of an homogeneous element $x$.
$A$ is said to be {\it commutative} if 
$$
xy=(-1)^{p(x)p(y)}yx
$$
and its category will be denoted by $\salg$. 

\medskip

The concept of an affine supervariety or more generally an
affine superscheme can be defined very effectively 
through its functor of points. 

\medskip

\begin{definition} \label{superscheme}
An {\it affine superscheme} is a
representable functor:
$$
\begin{array}{cccc}
\bX: & \salg & \lra    & \sets \\ \\
   & A       & \mapsto & X(A)=\Hom(k[X], A) 
\end{array}
$$

\end{definition}

\medskip

From this definition one can see that 
the category $\salg$ plays a role in
algebraic supergeometry similar to the category of
commutative algebras for the ordinary (i.e. non super)
algebraic geometry. In particular it is possible to
show that there is an equivalence of categories between the categories
of affine superschemes and commutative superalgebras. 
(For more details see \cite{fl}).

\medskip

\begin{examples} \label{svs}

\medskip
\
1. {\it Affine superspace}. Let $V=V_0 \oplus V_1$ 
be a finite dimensional super vector space.
Define the following functor:
$$
\bV: \salg \lra \sets , \qquad 
\bV(A)=(A \otimes V)_0= A_0 \otimes V_0 \oplus A_1 \otimes V_1
$$
This functor is representable and it is represented by:
$$
k[V]=Sym(V_0)\otimes \wedge(V_1)
$$
where $Sym(V_0)$ is the polynomial algebra over the vector space $V_0$ and
$\wedge(V_1)$ the exterior algebra over the vector space $V_1$.
Let's see this more in detail.

\medskip

If we choose a graded basis for $V$, $e_1 \dots e_m, \ep_1 \dots \ep_n$,
with $e_i$ even and $\ep_j$ odd,
then 
$$
k[V]=k[x_1 \dots x_m,\xi_1 \dots \xi_n],
$$ 
where
the latin letters denote commuting indeterminates, while the greek ones 
anticommuting indeterminates i.e. $\xi_i \xi_j=-\xi_j\xi_i$.
In this case $V$ is commonly denoted with $k^{m|n}$
and $m|n$ is called the {\it superdimension} of $V$.
We also will call $\bV$ as the {\it functor of points} of the
super vector space $V$.

Observe that:
$$
\begin{array}{c}
\bV(A)=\left\{
(a_1 \dots a_m, \al_1 \dots \al_n)
\quad | \quad a_i \in A_0, \quad \al_j \in A_1 \right\}  = \\ \\
\Hom(k[V],A)=\{ \phi:k[V] \lra k \quad | \quad \phi(x_i)=a_i, \quad
\phi(\xi_j)=\al_j\}
\end{array} 
$$

Hence $\bV(A)= A_0 \otimes k^m \oplus A_1 \otimes k^n$.

\medskip


\medskip

2. {\it Tensor superspace}.
We define the vector space of $r$-tensors as:
$$
T^r(V)=_{\defi} \underbrace{V \otimes V \dots \otimes
V}_{\hbox{$r$ times}}
$$
$T^r(V)$ is a super vector space, the parity of a
monomial element is defined as
$p(v_1 \otimes \dots \otimes v_r)=
p(v_1)+ \dots +p(v_r)$.
$T^r(V)$ is also a supervariety functor:
$$
{\bf T^r(V)}(A)=\bV(A) \otimes_A \dots \otimes_A \bV(A)
$$
We define the superspace of tensors $T(V)$ as:
$$
T(V)=\bigoplus_{r \geq 0} T^r(V)
$$
and denote with ${\bf \bT(V)}$ its functor of points.
\medskip

3. {\it Supermatrices}. 
Given a finite dimensional super vector space $V$, the endomorphisms
$\End(V)$ over $V$ is itself a supervector space: 
$\End(V)=\End(V)_0 \oplus \End(V)_1$, where $\End(V)_0$ are the endomorphisms
preserving parity, while $\End(V)_1$ are those reversing parity.
 
Hence we can define the following functor:
$$
\bEnd(V): \salg \lra \sets, \qquad
\bEnd(V)(A)=(A \otimes \End(V))_0
$$

This functor is representable (see (1)). Choosing a graded basis for $V$,
$V=k^{m|n}$, the functor is represented by
$k[x_{ij},y_{il},\xi_{kj},\eta_{kl}]$
where $1 \leq i,j \leq m$, $m+1 \leq k,l \leq m+n$.

In this case:
$$
\bEnd(V)(A)=\left\{ 
\begin{pmatrix} 
a_{m\times m} & \beta_{m\times n} \\ 
\gamma_{n \times m} & d_{n\times n} 
\end{pmatrix}
\right\}
$$
where $a$, $d$ and $\beta$, $\gamma$ are block matrices with
respectively even and odd entries. 

\end{examples}
 
\medskip
 
\begin{definition} An {\it affine supergroup $G$} is a group valued affine
superscheme, i.e. it is a representable functor:
$$
\begin{array}{cccc}
\bG: & \salg & \lra & \groups \\
& A & \mapsto & \bGL(V)(A)
\end{array}
$$

\medskip
It is simple to verify that the superalgebra representing the
supergroup $\bG$ has an Hopf superalgebra structure. More is true:
Given a supervariety $\bG$,
$\bG$ is a supergroup if and only if the algebra representing
it $k[G]$ is an
Hopf superalgebra.
\end{definition}

\medskip

Let $V$ be a finite dimensional super vector space.
We are interested in the {\it general linear supergroup} $\bGL(V)$. 

\medskip


\begin{definition}   
We define {\it general linear supergroup} the functor
$$
\begin{array}{cccc}
\GL(V): & \salg & \lra & \sets \\ \\
& A & \mapsto & \bGL(V)(A)
\end{array}
$$ 
where $\bGL(V)(A)$ is the set of automorphisms of the $A$-supermodule
$A \otimes V$, $A \in \salg$.
More explicitly if $V=k^{m|n}$, the functor $\bGL(V)$ commonly denoted
$\bGL(m|n)$ is defined as the set of automorphisms of 
$A^{m|n}=_{\defi} A \otimes k^{m|n}$ and is given by:

$$
\bGL(m|n)(A)=\left\{ 
\begin{pmatrix} 
a_{m\times m} & \beta_{m\times n} \\ 
\gamma_{n \times m} & d_{n\times n} 
\end{pmatrix} \quad |\quad
a, d \quad \hbox{invertible} 
\right\}
$$
where $a$, $d$ and $\beta$, $\gamma$ are block matrices with
respectively even and odd entries. 

\medskip

This functor is representable and it is represented by 
the Hopf algebra (see \cite{fi1}):
$$
\begin{matrix}
k[x_{ij}, y_{\alpha \beta},
\xi_{i\beta},\eta_{\alpha
j},z,w]/\bigr((w\det(x)-1,z\det(y)-1\bigl),
\cr\cr i,j=1,\dots m \qquad 
\alpha,\beta=1,\dots n.
\end{matrix}
$$
\end{definition}

We now would like to introduce the notion of Lie superalgebra using
the functorial language. We then see it is equivalent to
the more standard definitions (see \cite{ka} for example).

\medskip

\begin{definition} \label{superlie}
Let $\fg$ be a finite dimensional supervector space. 
The functor (see Example \ref{svs} (1)):
$$
\bg: \salg \lra \sets, \qquad 
\bg(A)=(A \otimes \bg)_0 
$$
is said to be a {\it Lie superalgebra} if it is
Lie algebra valued, i.e. for each $A$ there exists a 
linear map:
$$
[\;,\;]_A:\bg(A) \times \bg(A) \lra \bg(A)
$$ 
satisfying the 
antisymmetric property and the
Jacobi identity.
\end{definition}

\medskip

Notice that in the same way as the supergroup functor is group valued, the
Lie superalgebra functor is Lie algebra valued, i. e. it has
values in a {\sl classical category}. The super nature
of these functors arises from the different starting category,
namely $\salg$, which allows superalgebras as representing objects.

\medskip

The usual notion
of Lie superalgebra, as defined for example by Kac in \cite{ka}
is equivalent to this functorial definition.
Let's recall this definition and see the equivalence with the Definition
\ref{superlie} more in detail.

\medskip

\begin{definition}
Let $\fg$ be a super vector space. We say that
a bilinear map 
$$
[,]: \fg \times \fg \lra \fg
$$
is a {\it superbracket} if  $\forall x,y,z \in \fg$:
$$
[x,y] = (-1)^{p(x)p(y)}[y,x] 
$$
$$
[x,[y,z]] + (-1)^{p(x)p(y)+p(x)p(z)}[y,[z,x]]+ 
 (-1)^{p(x)p(z)+p(y)p(z)} [z,[x,y]]=0
$$

$(\fg,[,])$,  is what in the literature is commonly defined as 
{\sl Lie superalgebra}.

\end{definition}

\medskip

\begin{observation} 
The two concepts of Lie superalgebra $\bg$ in the
functorial setting and superbracket
on a supervector space $(\fg,[,])$ are equivalent.


\medskip

In fact if we have a Lie superalgebra $\bg$
there is always a superspace
$\fg$ associated to it together with a superbracket.
The superbracket on $\fg$ is given following the {\it even rules}. 
(For a complete treatment of even rules see pg 57 \cite{dm}).
\medskip
Given $v, w \in \fg$, we have that since the Lie bracket on $\bg(A)$ is
$A$-linear:
$$
[a \otimes v, b \otimes w]=ab \otimes z \in (A \otimes \fg)_0 = \fg(A) 
$$
Hence we can define the bracket $\{v,w\}$ as the element of $\fg$ such that:
$z=(-1)^{p(a)p(w)}\{v,w\}$ i. e. satisfying the relation:
$$
[a \otimes b, b \otimes w]= (-1)^{p(b)p(v)}ab \otimes \{v,w\} 
$$
 
We need to check it si a superbracket. Let's see
for example the antisymmetry property.
Observe first that if $a \otimes v \in (\fg \otimes A)_0$ 
must be $p(v)=p(a)$, since $(A \otimes \fg)_0=$
$A_0 \otimes \fg_0  \oplus A_1 \otimes \fg_1$.
So we can write:
$$
[a \otimes v, b \otimes w]=(-1)^{p(b)p(v)}ab \otimes \{v,w\} =
(-1)^{p(v)p(w)}ab \otimes \{v,w\}  
$$
On the other hand:

\begin{eqnarray*}
[b \otimes w, a \otimes v] &=(-1)^{p(a)p(w)}ba \otimes \{w,v\} = \\ \\
&=(-1)^{p(a)p(w)+p(a)p(b)}ab \otimes \{w,v\} = \\ \\
&=(-1)^{2p(w)p(v)} ab \otimes \{w,v\} =ab \otimes \{w,v\}.
\end{eqnarray*}

Comparing the two expression we get the antisymmetry of the superbracket.
For the super Jacobi identity the calculation is the same.

\medskip

Vice versa if $(\fg, \{,\})$ is a super vector space with a superbracket,
we immediately can define its functor of points $\bg$. $\bg$ is a Lie
superalgebra because we have a bracket on $\bg(A)$ defined as
$$
[a \otimes v, b \otimes w]=(-1)^{p(b)p(v)} ab \otimes \{v,w\}
$$
The previous calculation worked backwards proves that $[,]$ is Lie
bracket.

\end{observation}

\medskip

With an abuse of language we will call Lie superalgebra  both the supervector
space $\fg$ with a superbracket $[,]$ and the functor $\bg$  as
defined in \ref{superlie}.

\medskip
\begin{observation}
In \cite{fl} is given the notion of a Lie super algebra associated
to an affine supergroup. In this work we show
that the Lie superalgebra associated to $\bGL(m|n)$ is $\bEnd(k^{m|n})$.
We will denote $\bEnd(k^{m|n})$ with
$\fgl(m|n)$ as supervector space
and with $\bgl(m|n)$ as its functor of points.
The purpose of this paper does not allow for a full description of such
correspondence, 
all the details and the proofs can be found in 
\cite{fl}. 
\end{observation}
\bigskip
\section{Summary and observations on results by Berele and Regev} \label{berele}

\medskip

In this section we want to review  some of the  results in \cite{br,dj}.
We wish to describe the correspondence
between tensorial representations of the superalgebra
$\fgl(m|n)$ and representations of the permutation
group. 
This correspondence is obtained using
the double centralizer theorem.
(Note: in \cite{br} $\fgl(m|n)$ is
denoted by $\mathfrak{pl}$).

\medskip

Let $V=k^{m|n}$ 
and let $T(V)=\bigoplus_{r \geq 0} T^r(V)$ be the tensor superspace
(see Example \ref{svs} (2)).

\medskip

We want to define on $T^r(V)$ two actions: one by $S_r$ the
permutation group and the other by the Lie superalgebra $\fgl(m|n)$.

\medskip

Let $\sigma=(i,j) \in S_r$
and let $\{v_i\}_{1\leq i \leq m+n}$ be a 
basis of $V$ ($v_1 \dots v_m$ even elements and
$v_{m+1} \dots v_{m+n}$ odd ones. Let's define:
$$
(v_1 \otimes \dots \otimes v_r) \cdot \sigma=_{\defi}
\epsilon v_{\s(1)} \otimes \dots \otimes v_{\s(r)}
$$
where $\epsilon=-1$ when $v_i$ and $v_j$ are both odd
and $\epsilon=1$ otherwise. 
This defines a representation $\tau_r$ of $S_r$ in $T^r(V)$. The
proof of this fact can be found in \cite{br} pg 122-123.

\medskip

Consider now the action $\theta_r$ of the Lie superalgebra $\fgl(m|n)$ on
$T^r(V)$ given by derivations:
$$
\begin{matrix}
\theta_n(X) 
(v_1 \otimes \dots \otimes v_r)=_{\defi} 
\sum_i (-1)^{s(X,i)} v_1 \otimes X(v_i) \otimes \dots
\otimes v_r  \\ \\ 
\qquad g \in \fgl(V)(A), \quad v_i \in V(A), \quad A \in \salg 
\end{matrix}
$$ 
with $s(X,i)=p(X)o(i)$ where $o(i)$ denotes the number of odd
elements among $v_1 \dots v_i$.

One can see that this is a Lie superalgebra
action i.e. it preserves the superbracket
and that it extends to an action $\theta$ of
$\fgl(m|n)$ on $T(V)$ (this is proved in \cite{br} 4.7).

\medskip

In \cite{br} Theorem 4.14 and Remark 4.15 is proved the important
double centralizer theorem:

\begin{theorem}
The algebras $\tau(S_r)$ and $\theta(\fgl(m|n))$ are
each the centralizer of the other in $\End(T^r(V))$.
\end{theorem}

\medskip

This result establishes a one to one correspondence
between irreducible tensor representations of $S_r$ occurring in
$\tau_r$ and those of $\fgl(m|n)$ occurring in $\theta_r$.

\medskip

These representations are parametrized by partitions $\lambda$
of the integer $r$. In \cite{br} \S 3 and \S 4 is worked out
completely the structure of irreducible tensor representations
of $\fgl(m|n)$ arising in this way. We are interested in their
dimensions.

\medskip

\begin{definition} \label{semistandard}
Let $t_1 < \dots < t_m < u_1 < \dots < u_n$ be integers
and $\lambda$ a partition of $r$ corresponding to a diagram $D_\la$.
A filling $T_\la$ of $D_\la$ is a {\it semistandard or superstandard tableau} if

1. The part of $T_\la$ filled with the $t$'s is a tableaux.

2. The $t$'s are non decreasing in rows and strictly increasing in columns.

3. The $u$'s are  non decreasing in columns and strictly increasing in rows.  
\end{definition}

As an example that will turn out to be important later let's look at
$m=n=1$, $t_1=1$, $u_1=2$ and $r=2$. We can have only two partitions: 
$\la=(2)$, $\la=(1,1)$. Each partition admits two fillings:
$$
\begin{array}{ccccc}
\la=(2)   & & \begin{array}{cc} 1 & 1 \end{array} & & 
\begin{array}{cc} 1 & 2 \end{array} \\ \\
\la=(1,1) & & \begin{array}{c} 1 \\ 2 \end{array} & & 
\begin{array}{c} 2 \\ 2 \end{array}
\end{array}
$$

By Theorem 3.17, 3.18 and 4.17 in \cite{br} we have the following:

\medskip
\begin{theorem} \label{bereleregev}
The irreducible representations of $\fgl(m|n)$ occurring in $\theta_r$ 
are parametrized by partitions of $\lambda$ of the integer $r$.
The irreducible representations associated to the shape $\lambda$
has dimension equal to the number of semistandard tableaux of
shape $\lambda$.
\end{theorem}

\medskip

\begin{observation} \label{trace}
This theorem tells us immediately that we have no one dimensional
representations of $\fgl(m|n)$ occurring in $\theta_r$, if $n >0$.
In fact one can generalize the Example \ref{semistandard} to show that
since the odd variables allow repetitions on rows, we  
always have more than one filling for each shape. However there
exists a polynomial representation of $\fgl(m|n)$ of dimension one,
namely the supertrace (\cite{be} pg. 100):
$$
\begin{array}{ccc}
\fgl(m|n) & \lra & k \cong \End(k) \\
A=
\begin{pmatrix} 
X & Y \\ Z & W 
\end{pmatrix}
& \mapsto & str(A)=_{\defi}tr(X)-tr(W)
\end{array}
$$
This shows that the tensor representations described in \cite{br}
do not exhaust all polynomial representations of $\fgl(m|n)$, for $n>0$.
\end{observation}

%
%
%
%
%

\section{Tensor representations of the general linear supergroup} 

\label{supergroup}

\medskip

Let's start by introducing the notion of supergroup
and of Lie super algebra representation
from a functorial point of view.

\medskip

\begin{definition}
Given an affine algebraic supergroup $\bG$ we say that 
$\bG$ acts on a super vector space $W$, if we have a 
natural transformation:
$$
\rr: \bG \lra \bEnd(W)
$$
In other words, if we have for a fixed $A \in \salg$ a
functorial morphism $\rr_{A}: \bG(A) \lra \bEnd(W)(A)$. 
If $W \cong k^{m|n}$ we can identify $\rr_{A}(g)$ with
a matrix in $\bEnd(W)(A)$ (see Example \ref{svs} (3)).
\end{definition}

\medskip

Let $V$ be a finite dimensional super vector space. Define: 
$$
\begin{matrix}
\rho_r:\bGL(V) \lra \bEnd(T^r(V))
\\ \\
\rho_{r,A}(g) v_1 \otimes \dots \otimes v_n=_{\defi} g(v_1) \otimes \dots
\otimes g(v_n), \\ \\
g \in \bGL(V)(A), \quad v_i \in \bV(A), \quad A \in \salg 
\end{matrix}
$$
This is an action of $\bGL(V)$ on $T^r(V)$, that can be
easily extended to the whole $T(V)$.

\medskip

We now introduce the concept of a Lie superalgebra representation
using the functorial language.

\medskip

\begin{definition} 
Given a Lie superalgebra $\bg$ we say that 
$\bg$ acts on a super vector space $W$, if we have a 
natural transformation:
$$
t: \bg \lra \bEnd(W)
$$
preserving the Lie bracket, that is 
for a fixed $A \in \salg$, we have a Lie algebra 
morphism $t_A: \bg(A) \lra \End(W)(A)$.
It is easy to verify that this is equivalent to ask
that we have a morphism of Lie superalgebras:
$$
T: \fg \lra \End(W)
$$
i.e. a super vector space morphism preserving the superbracket. 
This agrees with the definition of Lie superalgebra representation
in \cite{br}, which we also recalled in \S 3.
\end{definition}

\medskip

We are interested in the action $\theta_r$ of $\fgl(V)$, the Lie superalgebra of
$\GL(V)$ on $T^r(V)$ introduced in Section \ref{berele}.


\medskip

Let's assume from now on $V=k^{m|n}$.
Denote with $\{\re_{ij}\}$ the graded canonical basis for
the supervector space $\fgl(m|n)$, with
$p(\re_{ij})=p(i)+p(j)$.

\medskip

\begin{definition}
Consider the following functor $\bE_{ij}: \salg \lra \sets$,
$1 \leq i \neq j \leq m+n$:
$$
\bE_{ij}(A)=\{I+x\re_{ij} | x \in A_k, k=p(\re_{ij})\} 
$$
This is an affine supergroup functor represented by
$k^{1|0}$ if $p(i)+p(j)$ is even, by $k^{0|1}$ if it is odd.
We call $\bE_{ij}$ a {\it one parameter subgroup functor}.
\end{definition}

\medskip

Consider also the functor $\bH_i: \salg \lra \sets$:
$$
\bH_{i}(A)=\{I+(x-1)\re_{ii} | x \in A_0 \setminus \{0\} \} 
$$
This is also an affine supergroup functor represented by
$(k^{1|0})^\times$, the multiplicative group of the ground
field $k$.

\medskip 


\medskip

\begin{theorem}
1. The affine supergroup functor $\rGL(m|n)$ is generated by
the subgroup functors $\{\bE_{ij}, \bH_i\}$, that is
the group $GL(m|n)(A)$ is generated by $\{\bE_{ij}(A), \bH_i(A)\}$ for all
$A \in \salg$.

2. The Lie superalgebra $\bgl(m|n)$ is generated by the
functors $\bfe_{ij}$ where:
$$
\bfe_{ij}=\{ a \otimes \re_{ij} |p(a)=p(i)+p(j) \}
$$
\end{theorem}

\Proof. (2) is immediate.
For (1) it is enough to prove $\bE_{ij}(A)$ generate the following
(see \cite{vsv} pg. 117): 
$$
\begin{pmatrix} 
A & 0 \\ 0 & D 
\end{pmatrix}
\qquad
\begin{pmatrix} 
I & B \\ 0 & I 
\end{pmatrix}
\qquad
\begin{pmatrix} 
I & 0 \\ C & I 
\end{pmatrix}
$$
The fact they generate the first type of matrices comes from
the classical theory. The fact they generate the other types
is immediate. \qed

\medskip

Let $\rGL(m)$ and $\rGL(n)$ denote the general linear
group of the ordinary vector spaces $V_0=k^m$ and $V_1=k^n$. 
Consider now the action of the (non super)
group $\GL(m) \times \GL(n)$ on the ordinary vector
space $V=V_0 \oplus V_1$ and also the action of its Lie
algebra $\fgl(m) \times \fgl(n)$ on the same space.
We can build the diagonal action $\rho^0$ of 
$\GL(m) \times \GL(n)$ on the space of tensors
$T(V)$ (again $V$ is viewed disregarding the grading)
and also the usual action $\theta^0$ by derivation of 
$\fgl(m) \times \fgl(n)$ on the same space.

\medskip

\begin{lemma} \label{classical}
$$
<\rho^0(\rGL(m) \times \rGL(n))>=
<\theta^0(\fgl(m) \times \fgl(n))>
$$
where $<S>$ denotes the subalgebra generated by the set $S$ inside
$\End(V)$ the endomorphism of the ordinary vector space $V=k^{m+n}$.
\end{lemma}

\Proof. This is a consequence of a classical result,
see for example \cite{jk} 8.2. \qed

\medskip


\begin{theorem} 
$$
<\rho_{r,A}(\GL(m|n)(A))>_A=<\theta_{r,A}(\fgl(m|n)(A))>_A \qquad A \in \salg.
$$
where $<S>_A$ denotes the subalgebra generated by the set $S$ inside
$\bEnd(V)(A)$.
\end{theorem}

\Proof. 
Since $\GL(m|n)$ is generated by $\{\bE_{ij}, \bH_i\}$
and $\fgl(m|n)$ is generated by $\{ \bfe_{ij}\}$ it is enough
to show that 
$$
\rho_{r,A}(\bE_{ij}(A)) \in \theta_{r,A}(\bgl(V)(A)), \qquad
\theta_{r,A}(\bfe_{ij}(A)) \in \rho_{r,A}(\bGL(V)(A))
$$
When $p(i)+p(j)$ is even this is an easy consequence of
Lemma \ref{classical}. 

Now the case when $p(i)+p(j)$ is odd. 
Let $D_{ij}$ be the derivation corresponding to 
the elementary matrix $\re_{ij}$. So we have
that $\bfe_{ij}(A)=\alpha_1\otimes D_{ij}$, $\alpha \in A_1$.
We claim that
$$
\theta_{r,A}(\bfe_{ij}(A))=1(A)-\rho_{r,A}(\bE_{ij})(A)
$$
This is a calculation.
\qed

\medskip
\begin{corollary} \label{GLrep}
There is a one to one correspondence between the irreducible representations
of $S_r$ and the irreducible representations of $\bGL(m|n)$ occurring
in $\rho_r$.
\end{corollary}

\medskip
\begin{observation}
By Corollary \ref{GLrep} and Theorem 
\ref{bereleregev} we have that also the irreducible
representations occuring in $\rho_r$ of $\bGL(m|n)$ are parametrized
by partitions of the integer $r$. However by Observation \ref{trace}
we have that there is no one dimensional irreducible representation
hence also for $GL(m|n)$ we miss an important representation, namely
the Berezinian:
$$
\begin{array}{ccc}
\bGL(m|n)(A) & \lra & A \cong \bEnd(k)(A) \\
\begin{pmatrix} 
X & Y \\ Z & W 
\end{pmatrix}
& \mapsto & det(W)^{-1}det(X-YW^{-1}Z)
\end{array}
$$
This shows that the tensor representations of $\bGL(m|n)$
so not exhaust all polynomial representations of $\bGL(m|n)$, for $n>0$.

\medskip

The Berezinian representation has been described by Deligne and
Morgan in \cite{dm} pg 60, in a natural way as
an action of $\bGL(V)$ on $\rExt$ group $\rExt_{\Sym^*(V^*)}^m
(A, \Sym^*(V^*))$.
$\rExt$ plays the same role as the antisymmetric tensors in this
super setting. It would be interesting to see if there are other
objects that give raise to representations which are not among those 
previously described.
\end{observation}


\bigskip

\end{document}